\documentclass[12pt]{amsart}
\usepackage{amsxtra}
\usepackage{amssymb, amsmath}
\usepackage{amsthm}
\usepackage {hyperref}

\newtheorem{theorem}{Theorem}

\newtheorem{corollary}{Corollary}

\thanks{This research is partially supported by the Natural Sciences and
Engineering Research Council's Discovery Grant to Pranesh Kumar}

\keywords{Relative information of type s; Relative JS --
divergence of type s; Relative AG -- divergence of type s;
Csisz\'{a}r $f$--divergence; Information inequalities.}
\subjclass[2000]{94A17; 26D15}

\begin{document}
\title[Information Measures]{On Unified Generalizations of Relative
Jensen--Shannon and Arithmetic--Geometric Divergence Measures, and
Their Properties}

\author{Pranesh Kumar}
\address{Pranesh Kumar\\
Mathematics Department\\
College of Science and Management\\
University of Northern British Columbia\\
Prince George BC V2N4Z9, Canada.} \email{kumarp@unbc.ca}
\urladdr{http://web.unbc.ca/$\sim$kumarp}

\author{Inder Jeet Taneja}
\address{Inder Jeet Taneja\\
Departamento de Matem\'{a}tica\\
Universidade Federal de Santa Catarina\\
88.040-900 Florian\'{o}polis, SC, Brazil}
\email{taneja@mtm.ufsc.br} \urladdr{http://www.mtm.ufsc.br/$\sim
$taneja}

\begin{abstract}
In this paper we shall consider one parametric generalization of
some non-symmetric divergence measures. The \textit{non-symmetric
divergence measures} are such as: Kullback-Leibler
\textit{relative information}, $\chi ^2-$\textit{divergence},
\textit{relative J -- divergence}, \textit{relative Jensen --
Shannon divergence} and \textit{relative Arithmetic -- Geometric
divergence}. All the generalizations considered can be written as
particular cases of Csisz\'{a}r's \textit{f-divergence}. By
putting some conditions on the probability distribution, the aim
here is to develop bounds on these measures and their parametric
generalizations.
\end{abstract}

\maketitle

\section{Introduction}

Let
\[
\Gamma _n = \left\{ {P = (p_1 ,p_2 ,...,p_n )\left| {p_i >
0,\sum\limits_{i = 1}^n {p_i = 1} } \right.} \right\}, \,\, n
\geqslant 2,
\]

\noindent be the set of all complete finite discrete probability
distributions. There exist many information and divergence
measures in the literature on information theory and statistics.
Some of them are symmetric with respect to probability
distributions, while others are not. Here, in this paper, we shall
work only with non-symmetric measures. Through out the paper it is
understood that the probability distributions $P,Q \in \Gamma _n
$.

\subsection{Non-Symmetric Measures}

Here we shall give some non-symmetric measures of information. The
well known among them are $\chi ^2 - $\textit{divergence} and
Kullback-Leibler \textit{relative information}.\\

$\bullet$ $\chi ^2 - $\textbf{Divergence} (Pearson \cite{pea})
\begin{equation}
\label{eq1} \chi ^2(P\vert \vert Q) = \sum\limits_{i = 1}^n
{\frac{(p_i - q_i )^2}{q_i }} = \sum\limits_{i = 1}^n {\frac{p_i^2
}{q_i } - 1}.
\end{equation}

$\bullet$ \textbf{Relative Information} (Kullback and Leiber
\cite{kul})
\begin{equation}
\label{eq3} K(P\vert \vert Q) = \sum\limits_{i = 1}^n {p_i \ln
(\frac{p_i }{q_i })}.
\end{equation}

$\bullet$ \textbf{Relative J--Divergence} (Dragomir et al.
\cite{dgp1})
\begin{equation}
\label{eq5} D(P\vert \vert Q) = \sum\limits_{i = 1}^n {(p_i - q_i
)\ln \left( {\frac{p_i + q_i }{2q_i }} \right)}.
\end{equation}

$\bullet$ \textbf{Relative JS--Divergence} (Sibson \cite{sib})
\begin{equation}
\label{eq7} F(P\vert \vert Q) = \sum\limits_{i = 1}^n {p_i \ln
\left( {\frac{2p_i }{p_i + q_i }} \right)}.
\end{equation}

$\bullet$ \textbf{Relative AG--Divergence} (Taneja \cite{tan6},
\cite{tan7})
\begin{equation}
\label{eq9} G(P\vert \vert Q) = \sum\limits_{i = 1}^n {\left(
{\frac{p_i + q_i }{2}} \right)\ln \left( {\frac{p_i + q_i }{2p_i
}} \right)}.
\end{equation}

The symmetric versions of the above measures are given by
\begin{equation}
\label{eq11} \Psi (P\vert \vert Q) = \chi ^2(P\vert \vert Q) +
\chi ^2(Q\vert \vert P),
\end{equation}
\begin{align}
\label{eq12} J(P\vert \vert Q) & = K(P\vert \vert Q) + K(Q\vert
\vert P)\\
&  = D(P\vert \vert Q) + D(Q\vert \vert P),\notag
\end{align}
\begin{equation}
\label{eq13} I(P\vert \vert Q) = \frac{1}{2}\left[ {F(P\vert \vert
Q) + F(Q\vert \vert P)} \right]
\end{equation}

\noindent and
\begin{equation}
\label{eq14} T(P\vert \vert Q) = \frac{1}{2}\left[ {G(P\vert \vert
Q) + G(Q\vert \vert P)} \right].
\end{equation}

After simplification, we can write
\begin{equation}
\label{eq15} J(P\vert \vert Q) = 4\left[ {I(P\vert \vert Q) +
T(P\vert \vert Q)} \right].
\end{equation}

Dragomir et al. \cite{dsb} studied the measures (\ref{eq11}). We
call it \cite{tan5} by \textit{symmetric chi-square divergence}.
The measure (\ref{eq11}) is well known Jeffreys-Kullback-Leiber
\cite{jef}, \cite{kul} \textit{J-divergence}. The measure
(\ref{eq13}) is \textit{information radius} studied by Sibson
\cite{sib}. It is also known by \textit{Jensen Shannon divergence
measure} (ref. Burbea and Rao \cite{bur}). The measure
(\ref{eq14}) is new in the literature and is studied for the first
time by Taneja \cite{tan2} and is called \textit{arithmetic and
geometric mean divergence measure}. Lin \cite{lin} studied some
interesting properties and applications of the measure \ref{eq7}.
More details on some of these divergence measures can be seen in
Taneja \cite{tan1}, \cite{tan2} and in on line book by Taneja
\cite{tan3}. An inequality among the measures
(\ref{eq11})-(\ref{eq14}) can be seen in Taneja \cite{tan8}.

In this paper our aim is to work with one parametric
generalization of non-symmetric divergence measures given by
(\ref{eq7}) and (\ref{eq9}). A similar kind of study of the
measures (\ref{eq3}) and (\ref{eq5}) with their one parametric
generalizations can be seen in Kumar and Taneja \cite{kut}.

It is important to note that defining a generalized divergence by
introducing a real parameter allows to unify many of the known
divergence measures studied individually and also yields a number
of new divergences. The properties and bounds established for this
family do hold good in particular for member divergences. It
provides a vast horizon of divergences for users to choose that
deems best for their applications. A few examples of new
divergence measures obtained from this generalized divergence are
cited in the next section 2.

\section{Generalized Non-Symmetric Divergence Measures}

A one parametric generalization of the measure ({eq1}) can be seen
in in Liese and Vajda \cite{liv}. We refer it here as \textit{relative information of type s}.\\

$\bullet$ \textbf{Relative Information of Type s}
\begin{equation}
\label{eq16} \Phi _s (P\vert \vert Q) = \begin{cases}
 {K_s (P\vert \vert Q) = \left[ {s(s - 1)} \right]^{ - 1}\left[
{\sum\limits_{i = 1}^n {p_i^s q_i^{1 - s} - 1} } \right],} & {s
\ne 0,1}
\\
 {K(Q\vert \vert P) = \sum\limits_{i = 1}^n {q_i \ln \left( {\frac{q_i }{p_i
}} \right)} ,} & {s = 0} \\
 {K(P\vert \vert Q) = \sum\limits_{i = 1}^n {p_i \ln \left( {\frac{p_i }{q_i
}} \right)} ,} & {s = 1} \\
\end{cases},
\end{equation}

\noindent for all $s \in \mathbb{R}$.

\bigskip
The measure (\ref{eq16}) admits the following particular cases:
\begin{itemize}
\item[(i)] $\Phi _{ - 1} (P\vert \vert Q) = \frac{1}{2}\chi
^2(Q\vert \vert P)$. \item[(ii)] $\Phi _0 (P\vert \vert Q) =
K(Q\vert \vert P)$. \item[(iii)] $\Phi _{1 / 2} (P\vert \vert Q) =
4\left[ {1 - B(P\vert \vert Q)} \right] = 4h(P\vert \vert Q)$.
\item[(iv)] $\Phi _1 (P\vert \vert Q) = K(P\vert \vert Q)$.
\item[(v)] $\Phi _2 (P\vert \vert Q) = \frac{1}{2}\chi ^2(P\vert
\vert Q)$.
\end{itemize}

Thus we observe that $\Phi _2 (P\vert \vert Q) = \Phi _{ - 1}
(Q\vert \vert P)$and $\Phi _1 (P\vert \vert Q) = \Phi _0 (Q\vert
\vert P)$.\\

The measures $B(P\vert \vert Q)$ and $h(P\vert \vert Q)$ appearing
in part (iii) are given by
\begin{equation}
\label{eq17} B(P\vert \vert Q) = \sum\limits_{i = 1}^n \sqrt {p_i
q_i }
\end{equation}

\noindent and
\begin{equation}
\label{eq18} h(P\vert \vert Q) = 1 - B(P\vert \vert Q) =
\frac{1}{2}\sum\limits_{i = 1}^n {(\sqrt {p_i } - \sqrt {q_i }
)^2}
\end{equation}

\noindent respectively.\\

The measure $B(P\vert \vert Q)$ is famous as Bhattacharyya
\cite{bha} \textit{coefficient} and the measure $h(P\vert \vert
Q)$is known as Hellinger \cite{hel} \textit{discrimination}.\\

Now we shall give one parametric generalization of the measures
given by (\ref{eq5}) and (\ref{eq7}). These generalizations are
based on the measure (\ref{eq16}).\\

$\bullet$ \textbf{Unified Relative AG and JS -- Divergence of Type
s}

\bigskip
Let us consider the following unified one parametric
generalizations of the measures (\ref{eq7}) and (\ref{eq9})
simultaneously.

\begin{equation}
\label{eq19} \Omega _s (P\vert \vert Q) = \begin{cases}
 {FG_s (P\vert \vert Q) = \left[ {s(s - 1)} \right]^{ - 1}\left[
{\sum\limits_{i = 1}^n {p_i \left( {\frac{p_i + q_i }{2p_i }}
\right)^s - 1}
} \right],} & {s \ne 0,1} \\
 {F(P\vert \vert Q) = \sum\limits_{i = 1}^n {p_i \ln \left( {\frac{2p_i
}{p_i + q_i }} \right)} ,} & {s = 0} \\
 {G(P\vert \vert Q) = \sum\limits_{i = 1}^n {\left( {\frac{p_i + q_i }{2}}
\right)\ln \left( {\frac{p_i + q_i }{2p_i }} \right)} ,} & {s = 1}
\\
\end{cases}.
\end{equation}

The measure (\ref{eq19}) admits the following particular cases:
\begin{itemize}
\item[(i)] $ \Omega _{ - 1} (P\vert \vert Q) = \frac{1}{4}\Delta
(P\vert \vert Q)$

\item[(ii)]$ \Omega _0 (P\vert \vert Q) = F(P\vert \vert Q).$

\item[(iii)]$ \Omega _{1 / 2} (P\vert \vert Q) = 4\left[ {1 -
B\left( {P\vert \vert \frac{P + Q}{2}} \right)} \right] =
4\,h\left( {P\vert \vert \frac{P + Q}{2}}\right).$

\item[(iv)]$ \Omega _1 (P\vert \vert Q) = G(P\vert \vert Q).$

\item[(v)]$ \Omega _2 (P\vert \vert Q) = \frac{1}{8}\chi ^2(Q\vert
\vert P).$
\end{itemize}

\bigskip
The expression $\Delta (P\vert \vert Q)$ appearing in part (i) is
the well known \textit{triangular discrimination}, and is given by
\begin{equation}
\label{eq21} \Delta (P\vert \vert Q) = \sum\limits_{i = 1}^n
{\frac{(p_i - q_i )^2}{p_i + q_i }} .
\end{equation}

The new divergences can be obtained from $\Omega _s (P\vert \vert
Q)$ by considering different choices of the real parameter $s$.
For example, $s = - \frac{1}{2}$ and $s = - 2$ in (\ref{eq19})
respectively result the following new divergence measures:
\begin{equation}
\Omega _{ - 1 / 2} (P\vert \vert Q) = \frac{4}{3}\left[
{\sum\limits_{i = 1}^n {p_i \sqrt {\frac{2p_i }{p_i + q_i }} } -
1} \right]
\end{equation}

\noindent and
\begin{equation}
\Omega _{ - 2} (P\vert \vert Q) = \frac{1}{6}\left[
{\sum\limits_{i = 1}^n {p_i \left( {\frac{2p_i }{p_i + q_i }}
\right)} ^2 - 1} \right].
\end{equation}

In this paper we shall study some properties of the
\textit{unified generalized measure} (\ref{eq19}). Some properties
and application of the this measure following the lines of Lin
\cite{lin} shall be dealt elsewhere. Some applications of this new
generalized measure (\ref{eq19}) towards, \textit{pattern
recognition}, \textit{statistics}, \textit{minimization problem},
etc. are also under study.

\section{Csisz\'{a}r's $f-$Divergence and its Properties}

Given a function$f:[0,\infty ) \to \mathbb{R}$, the
\textit{f-divergence} measure introduced by Csisz\'{a}r's
\cite{csi1} is given by
\begin{equation}
\label{eq23} C_f (P\vert \vert Q) = \sum\limits_{i = 1}^n {q_i
f\left( {\frac{p_i }{q_i }} \right)} ,
\end{equation}

\noindent for all $P,Q \in \Gamma _n $.

Here below are some theorems giving properties of the measure
(\ref{eq23}).

\begin{theorem} \label{the31} $($Csisz\'{a}r's \cite{csi1,csi2}$)$ If the function $f$
is convex and normalized, i.e., $f(1) = 0$, then the
\textit{f-divergence}, $C_f (P\vert \vert Q)$ is nonnegative and
convex in the pair of probability distribution $(P,Q) \in \Gamma
_n \times \Gamma _n $.
\end{theorem}

\begin{theorem} \label{the32} $($Dragomir \cite{dra1, dra2}$)$. Let $f:\mathbb{R}_ +
\to \mathbb{R}$ be differentiable convex and normalized i.e.,
$f(1) = 0$. Then
\begin{equation}
\label{eq24}
0 \leqslant C_f (P\vert \vert Q) \leqslant E_{C_f } (P\vert \vert Q)
\end{equation}

\noindent where
\begin{equation}
\label{eq25}
E_{C_f } (P\vert \vert Q) = \sum\limits_{i = 1}^n {(p_i - q_i )}
{f}'(\frac{p_i }{q_i }),
\end{equation}

\noindent for all $P,Q \in \Gamma _n $.

Let $P,Q \in \Gamma _n $ be such that there exists $r,R$ with $0 < r
\leqslant \frac{p_i }{q_i } \leqslant R < \infty $, $\forall i \in
\{1,2,...,n\}$, then
\begin{equation}
\label{eq26}
0 \leqslant C_f (P\vert \vert Q) \leqslant A_{C_f } (r,R),
\end{equation}

\noindent where
\begin{equation}
\label{eq27}
A_{C_f } (r,R) = \frac{1}{4}(R - r)\left[ {{f}'(R) - {f}'(r)} \right].
\end{equation}

Further, if we suppose that $0 < r \leqslant 1 \leqslant R <
\infty $, $r \ne R$, then
\begin{equation}
\label{eq28}
0 \leqslant C_f (P\vert \vert Q) \leqslant B_{C_f } (r,R),
\end{equation}

\noindent where
\begin{equation}
\label{eq29}
B_{C_f } (r,R) = \frac{(R - 1)f(r) + (1 - r)f(R)}{R - r}.
\end{equation}
\end{theorem}

Moreover, the following inequalities hold:
\begin{equation}
\label{eq30}
E_{C_f } (P\vert \vert Q) \leqslant A_{C_f } (r,R),
\end{equation}
\begin{equation}
\label{eq31}
B_{C_f } (r,R) \leqslant A_{C_f } (r,R)
\end{equation}

\noindent and
\begin{equation}
\label{eq32}
0 \leqslant B_{C_f } (r,R) - C_f (P\vert \vert Q) \leqslant A_{C_f } (r,R).
\end{equation}

The inequalities (\ref{eq30}) and (\ref{eq32}) can be seen in
Dragomir \cite{dra2}, while the inequality (\ref{eq31}) can be
proved easily.

\begin{theorem} \label{the33} Let $P,Q \in \Gamma _n $ be such that $0
< r \leqslant \frac{p_i }{q_i } \leqslant R < \infty $, $\forall i
\in \{1,2,...,n\}$, for some $r$ and $R$ with $0 < r < 1 < R <
\infty $. Let $f:\mathbb{R}_ + \to \mathbb{R}$ be differentiable
convex, normalized, of bounded variation, and second derivative is
monotonic with ${f}'''$ absolutely continuous on $[r,R]$ and
${f}''' \in L_\infty [r,R]$, then
\begin{align}
\label{eq33} & \left| {C_f (P\vert \vert Q) - \frac{1}{2}E_{C_f }
(P\vert \vert Q)} \right| \\
& \qquad \leqslant \min \left\{ {\frac{1}{8}k(f)\left[ {{f}''(R) -
{f}''(r)} \right]\chi ^2(P\vert \vert Q),} \right.\notag\\
& \qquad \qquad \qquad\left. {\frac{1}{12}\left\| {f}'''
\right\|_\infty \vert \chi \vert ^3(P\vert \vert Q),\,\left[
{{f}'(R) - {f}'(r)} \right]V(P\vert \vert Q)} \right\},\notag
\end{align}

\noindent and
\begin{align}
\label{eq34} & \left| {C_f (P\vert \vert Q) - E_{C_f }^\ast
(P\vert \vert Q)} \right|\\
& \qquad \leqslant \min \left\{ {\frac{1}{8}k(f)\left[ {{f}''(R) -
{f}''(r)} \right]\chi ^2(P\vert \vert Q),} \right.\notag\\
& \qquad \qquad \qquad\left. {\frac{1}{24}\left\| {f}'''
\right\|_\infty \vert \chi \vert ^3(P\vert \vert
Q),\,\frac{1}{2}\left[ {{f}'(R) - {f}'(r)} \right]V(P\vert \vert
Q)} \right\}, \notag
\end{align}

\noindent where
\begin{align}
\label{eq35} E_{C_f }^\ast (P\vert \vert Q) & = \sum\limits_{i =
1}^n {(p_i - q_i ){f}'\left( {\frac{p_i + q_i }{2q_i }} \right)}
,\\
\label{eq39} \vert \chi \vert ^3(P\vert \vert Q) & =
\sum\limits_{i = 1}^n {\frac{\vert p_i - q_i \vert ^3}{q_i^2 }},
\\
\label{eq43} V(P\vert \vert Q) & = \sum\limits_{i = 1}^n {\left|
{p_i - q_i } \right|}, \\
\label{eq36} \left\| {f}''' \right\|_\infty & = ess\mathop {\sup
}\limits_{x \in [r,R]} \left| {f}''' \right|,\\
\intertext{and} \label{eq35} k(f) & = \begin{cases}
 { - 1,} & {\mbox{if }{f}''\mbox{ is monotonically decreasing}} \\
 {1,} & {\mbox{if }{f}''\mbox{ is monotonically increasing}} \\
\end{cases}.
\end{align}
\end{theorem}

The above theorem is a combination of different papers due to
Dragomir et al. \cite{dgp1, dgp2, dgp3}.\\

The measures (\ref{eq1}), (\ref{eq39}) and (\ref{eq43}) are the
particular cases of Vajda \cite{vaj} $\vert \chi \vert ^m -
$\textit{divergence} given by
\begin{equation}
\label{eq44} \vert \chi \vert ^m(P\vert \vert Q) = \sum\limits_{i
= 1}^n {\frac{\vert p_i - q_i \vert ^m}{q_i^{m - 1} }} ,\mbox{ }m
\geqslant 1.
\end{equation}

The measure (\ref{eq44}) satisfy the following \cite{ced, dgp1}
properties
\begin{equation}
\label{eq45} \vert \chi \vert ^m(P\vert \vert Q) \leqslant
\frac{(1 - r)(R - 1)}{(R - r)}\left[ {(1 - r)^{m - 1} + (R - 1)^{m
- 1}} \right] \leqslant \left(\frac{R - r}{2}\right)^m, \,\, m
\geqslant 1.
\end{equation}

\noindent and
\begin{equation}
\label{eq46} \left( {\frac{1 - r^m}{1 - r}} \right)V(P\vert \vert
Q) \leqslant \vert \chi \vert ^m(P\vert \vert Q) \leqslant \left(
{\frac{R^m - 1}{R - 1}} \right)V(P\vert \vert Q), \,\, m \geqslant
1.
\end{equation}

For $m = 2$, $m = 3$ and $m = 1$ in (\ref{eq45}), we have
\begin{equation}
\label{eq54} \chi ^2(P\vert \vert Q) \leqslant (R - 1)(1 - r)
\leqslant \frac{(R - r)^2}{4},
\end{equation}
\begin{equation}
\label{eq55} \vert \chi \vert ^3(P\vert \vert Q) \leqslant
\frac{(R - 1)(1 - r)}{R - r}\left[ {(1 - r)^2 + (R - 1)^2} \right]
\leqslant \frac{1}{8}(R - r)^3
\end{equation}

\noindent and
\begin{equation}
\label{eq56} V(P\vert \vert Q) \leqslant \frac{2(R - 1)(1 - r)}{(R
- r)} \leqslant \frac{1}{2}(R - r).
\end{equation}
respectively.

In view of the last inequalities given in (\ref{eq54}),
(\ref{eq55}) and (\ref{eq56}) the bounds given in (\ref{eq33}) and
(\ref{eq34}) can be written in terms of $r,R$ as
\begin{align}
\label{eq57} & \left| {C_f (P\vert \vert Q) - \frac{1}{2}E_{C_f }
(P\vert \vert Q)} \right|\\
& \leqslant \frac{(R - r)^2}{4}\min \left\{ {\frac{1}{8}k(f)\left[
{{f}''(R) - {f}''(r)} \right],} \right. \, \left. {\frac{R -
r}{24}\left\| {f}''' \right\|_\infty ,\,\frac{2\left[ {{f}'(R) -
{f}'(r)} \right]}{R - r}} \right\}\notag
\end{align}

\noindent and
\begin{align}
\label{eq58} & \left| {C_f (P\vert \vert Q) - E_{C_f }^\ast
(P\vert \vert Q)} \right|\\
& \leqslant \frac{(R - r)^2}{4}\min \left\{ {\frac{1}{8}k(f)\left[
{{f}''(R) - {f}''(r)} \right],} \right. \,\left. {\frac{R -
r}{48}\left\| {f}''' \right\|_\infty ,\,\frac{{f}'(R) - {f}'(r)}{R
- r}} \right\},\notag
\end{align}

\noindent respectively.

\bigskip
From now onwards, unless otherwise specified, it is understood
that, if there are $r,R$, then $0 < r \leqslant \frac{p_i }{q_i }
\leqslant R < \infty $, $\forall i \in \{1,2,...,n\}$, with $0 < r
< 1 < R < \infty $, $P = (p_1 ,p_2 ,....,p_n ) \in \Gamma _n $ and
$Q = (q_1 ,q_2 ,....,q_n ) \in \Gamma _n $. Also, throughout the
paper we shall make use of the \textit{p-logarithmic power mean}
\cite{sto} given by
\begin{equation}
\label{eq59} L_p (a,b) = \begin{cases}
 {\left[ {\frac{b^{p + 1} - a^{p + 1}}{(p + 1)(b - a)}}
\right]^{\frac{1}{p}},} & {p \ne - 1,0} \\
 {\frac{b - a}{\ln b - \ln a},} & {p = - 1} \\
 {\frac{1}{e}\left[ {\frac{b^b}{a^a}} \right]^{\frac{1}{b - a}},} & {p = 0}
\\
\end{cases},
\end{equation}

\noindent for all $p \in \mathbb{R}$, $a \ne b$. In particular we
shall use the following notation
\begin{equation}
\label{eq60} L_p^p (a,b) = \begin{cases}
 {\frac{b^{p + 1} - a^{p + 1}}{(p + 1)(b - a)},} & {p \ne - 1} \\
 {\frac{\ln b - \ln a}{b - a}} & {p = - 1} \\
 1 & {p = 0} \\
\end{cases},
\end{equation}

\noindent for all $p \in \mathbb{R}$, $a \ne b$.

\section{Relative AG and JS -- Divergence of Type s}

In this section we shall consider the measures given by
(\ref{eq19}) and shall give some properties.

Let us consider
\begin{equation}
\label{eq61} \psi _s (x) = \begin{cases}
 {\left[ {s(s - 1)} \right]^{ - 1}\left[ {x\left( {\frac{x + 1}{2x}}
\right)^s - x - s\left( {\frac{1 - x}{2}} \right)} \right],} & {s \ne 0,1}
\\
 {\frac{1 - x}{2} - x\ln \left( {\frac{x + 1}{2x}} \right),} & {s = 0} \\
 {\frac{x - 1}{2} + \left( {\frac{x + 1}{2}} \right)\ln \left( {\frac{x +
1}{2x}} \right),} & {s = 1} \\
\end{cases},
\end{equation}

\noindent for all $x > 0$ in (\ref{eq23}), then $C_f (P\vert \vert
Q) = \Omega _s \left( {P\vert \vert Q} \right)$, where $\Omega _s
\left( {P\vert \vert Q} \right)$ is as given by (\ref{eq19}).

Moreover,
\begin{equation}
\label{eq62} \psi _s ^\prime (x) = \begin{cases}
 {(s - 1)^{ - 1}\left\{ {\frac{1}{s}\left[ {\left( {\frac{x + 1}{2x}}
\right)^s - 1} \right] + \frac{1}{2}\left[ {1 - \frac{1}{x}\left( {\frac{x +
1}{2x}} \right)^{s - 1}} \right]} \right\},} & {s \ne 0,1} \\
 {\frac{1 - x}{2(1 + x)} - \ln \left( {\frac{x + 1}{2x}} \right),} & {s = 0}
\\
 {\frac{1}{2}\left[ {1 - x^{ - 1} + \ln \left( {\frac{x + 1}{2x}} \right)}
\right],} & {s = 1} \\
\end{cases}
\end{equation}

\noindent and
\begin{equation}
\label{eq63} \psi _s ^{\prime \prime }(x) = \begin{cases}
 {\frac{1}{4x^3}\left( {\frac{x + 1}{2x}} \right)^{s - 2},} & {s \ne 0,1} \\
 {\frac{1}{x(x + 1)^2},} & {s = 0} \\
 {\frac{1}{2x^2(x + 1)},} & {s = 1} \\
\end{cases}
\end{equation}

Thus we have $\psi _s ^{\prime \prime }(x) > 0$ for all $x > 0$,
and hence, $\psi _s (x)$ is convex for all $x > 0$. Also, we have
$\psi _s (1) = 0$. In view of this we can say that\textit{
relative AG and JS -- divergence of type }$s$ is
\textit{nonnegative} and \textit{convex} in the pair of
probability distributions $(P,Q) \in \Gamma _n \times \Gamma _n $.

Based on Theorem \ref{the32}, we have the following theorem.

\begin{theorem} \label{the41} The following bounds on $\Omega _s
(P\vert \vert Q)$ hold:
\begin{equation}
\label{eq64}
\Omega _s (P\vert \vert Q) \leqslant E_{\Omega _s (P\vert \vert Q)} (P\vert
\vert Q) \leqslant A_{\Omega _s (P\vert \vert Q)} (r,R)
\end{equation}

\noindent and
\begin{equation}
\label{eq65}
\Omega _s (P\vert \vert Q) \leqslant B_{\Omega _s (P\vert \vert Q)} (r,R)
\leqslant A_{\Omega _s (P\vert \vert Q)} (r,R),
\end{equation}

\noindent where
\begin{align}
\label{eq66} E&_{\Omega _s (P\vert \vert Q)} (P\vert \vert Q) \\
& =
\begin{cases}
 {\left[ {s(s - 1)} \right]^{ - 1}\sum\limits_{i = 1}^n {\left( {\frac{p_i -
q_i }{p_i + q_i }} \right)\left( {\frac{p_i + q_i }{2p_i }} \right)^s\left[
{p_i + (1 - s)q_i } \right],} } & {s \ne 0,1} \\
 {D(Q\vert \vert P) - \frac{1}{2}\Delta (P\vert \vert Q),} & {s = 0} \\
 {\frac{1}{2}\left[ {\chi ^2(P\vert \vert Q) - D(Q\vert \vert P)} \right],}
& {s = 1} \\
\end{cases},\notag
\end{align}

\begin{align}
\label{eq67} A&_{\Omega _s (P\vert \vert Q)} (r,R)\\
& = \frac{(R - r)^2}{4rR}2^{ - s}\left\{ {L_{s - 1}^{s - 1} \left(
{\frac{r + 1}{r},\frac{R + 1}{R}} \right)} \right. \left. { - L_{s
- 2}^{s - 2} \left( {\frac{r + 1}{r},\frac{R + 1}{R}} \right)}
\right\}\notag
\end{align}

\noindent and
\begin{align}
\label{eq68} B&_ {\Omega _s (P\vert \vert Q)} (r,R)\\
& = \begin{cases}
 {\begin{array}{l}
 \frac{1}{2(s - 1)}L_{s - 1}^{s - 1} \left( {\frac{r + 1}{2r},\frac{R +
1}{2R}} \right) \\
 \mbox{ } + \frac{1}{s(s - 1)(R - r)}\left\{ {R\left[ {\left( {\frac{R +
1}{2R}} \right)^s - 1} \right] - r\left[ {\left( {\frac{r + 1}{2r}}
\right)^s - 1} \right]} \right\}, \\
 \end{array}} & {\begin{array}{l}
 \\
 \\
 s \ne 0,1 \\
 \end{array}} \\
 {\begin{array}{l}
 \\
 \frac{1}{(R - r)}\left[ {r\ln \left( {\frac{r + 1}{2r}} \right) - R\ln
\left( {\frac{R + 1}{2R}} \right)} \right] - \frac{1}{2}L_{ - 1}^{ - 1}
\left( {\frac{r + 1}{2r},\frac{R + 1}{2R}} \right), \\
 \\
 \end{array}} & {\begin{array}{l}
 \\
 s = 0 \\
 \\
 \end{array}} \\
 {\frac{rR - 1}{4rR}L_{ - 1}^{ - 1} \left( {\frac{r + 1}{2r},\frac{R +
1}{2R}} \right) + \frac{1}{2}\ln \left( {\frac{(R + 1)(r + 1)}{4rR}}
\right),} & {s = 1} \\
\end{cases}\notag
\end{align}
\end{theorem}

\begin{corollary} \label{cor41} The following inequality holds:
\begin{equation}
\label{eq69}
\frac{1}{2}\Delta (P\vert \vert Q) \leqslant D(Q\vert \vert P) \leqslant
\chi ^2(Q\vert \vert P).
\end{equation}
\end{corollary}

\begin{proof} It follows from (\ref{eq64}), where we take $s
= 0$ and $s = 1$ in (\ref{eq66}).
\end{proof}

Theorem \ref{the41} admits some particular cases. These are
summarized in the following two corollaries.

\begin{corollary} \label{cor42} The following bounds hold:
\begin{align}
\label{eq70} F(P\vert \vert Q) & \leqslant D(Q\vert
\vert P) - \frac{1}{2}\Delta (P\vert \vert Q)\\
& \leqslant \frac{(R - r)^2}{4rR}\left[ {L_{ - 1}^{ - 1} \left(
{\frac{r + 1}{r},\frac{R + 1}{R}} \right) - L_{ - 2}^{ - 2} \left(
{\frac{r + 1}{r},\frac{R + 1}{R}} \right)} \right]\notag
\end{align}

\noindent and
\begin{align}
\label{eq71} G(P\vert \vert Q) & \leqslant \frac{1}{2}\left[ {\chi
^2(Q\vert \vert P) - D(Q\vert \vert P)}
\right]\\
& \leqslant \frac{(R - r)^2}{8rR}\left[ {1 - L_{ - 1}^{ - 1}
\left( {\frac{r + 1}{r},\frac{R + 1}{R}} \right)} \right]\notag
\end{align}
\end{corollary}

\begin{proof} In inequalities (\ref{eq64}), take $s = 0$ and
$s = 1$ we get (\ref{eq70}) and (\ref{eq71}) respectively.
\end{proof}

For $s = - 1$ and $s = 2$ the results are obvious.

\begin{corollary} \label{cor43} The following bounds hold:
\begin{equation}
\label{eq72}
\Delta (P\vert \vert Q) \leqslant \frac{2(R - 1)(1 - r)}{(R + 1)(r + 1)}
\leqslant \frac{(R - r)^2(R + r + 2)}{(R + 1)^2(r + 1)^2},
\end{equation}
\begin{align}
\label{eq73} F(P\vert \vert Q) & \leqslant \frac{1}{R - r}\left[
{r\ln \left( {\frac{r + 1}{2r}} \right) - R\ln \left( {\frac{R +
1}{2R}} \right)} \right] - L_{ - 1}^{ - 1} \left( {\frac{r +
1}{r},\frac{R + 1}{R}} \right)\\
& \leqslant \frac{(R - r)^2}{4rR}\left[ {L_{ - 1}^{ - 1} \left(
{\frac{r + 1}{r},\frac{R + 1}{R}} \right) - L_{ - 2}^{ - 2} \left(
{\frac{r + 1}{r},\frac{R + 1}{R}} \right)} \right]\notag
\end{align}

\begin{align}
\label{eq74} G(P\vert \vert Q) &\leqslant \frac{rR - 1}{4rR}L_{ -
1}^{ - 1} \left( {\frac{r + 1}{2r},\frac{R + 1}{2R}} \right) +
\frac{1}{2}\ln \left( {\frac{(R + 1)(r + 1)}{4rR}} \right)\\
& \leqslant \frac{(R - r)^2}{8rR}\left[ {1 - L_{ - 1}^{ - 1}
\left( {\frac{r + 1}{r},\frac{R + 1}{R}} \right)} \right]\notag
\end{align}

\noindent and
\begin{equation}
\label{eq75}
\chi ^2(Q\vert \vert P) \leqslant \frac{(R - 1)(1 - r)}{rR} \leqslant
\frac{(R - r)^2(R + r)}{4r^2R^2}.
\end{equation}
\end{corollary}

\begin{proof} In inequalities (\ref{eq65}), take $s = - 1$,
$s = 0$, $s = 1$ and $s = 2$ we get respectively (\ref{eq72}),
(\ref{eq73}), (\ref{eq74}) and (\ref{eq75}).
\end{proof}

Based on Theorem \ref{the33}, we have the following result.

\begin{theorem} \label{the42} The following bounds hold:
\begin{align}
\label{eq76} & \left| {\Omega _s (P\vert \vert Q) -
\frac{1}{2}E_{\Omega _s } (P\vert \vert Q)}
\right|\\
& \leqslant \min \left\{ {\frac{1}{8}\delta _{\Omega _s }
(r,R)\chi ^2(P\vert \vert Q),\,\,\frac{1}{12}\left\| {\psi _s
^{\prime \prime \prime }} \right\|_\infty \left| \chi
\right|^3(P\vert \vert Q)}, \,\, \left[ {{f}'(R) - {f}'(r)}
\right] \right\}\notag\\
& \leqslant \frac{(R - r)^2}{32}\min \left\{ {\delta _{\Omega _s }
(r,R),\,\,\frac{R - r}{3}\left\| {\psi _s ^{\prime \prime \prime
}} \right\|_\infty },\,\,\frac{R - r}{2} \left[ {{f}'(R) -
{f}'(r)} \right] V(P\vert \vert Q) \right\}\notag
\end{align}

\noindent and
\begin{align}
\label{eq77} &\left| {\Omega _s (P\vert \vert Q) - E_{\Omega _s
}^\ast } (P\vert \vert
Q)\right|\\
& \leqslant \min \left\{ {\frac{1}{8}\delta _{\Omega _s }
(r,R)\chi ^2(P\vert \vert Q),\,\,\frac{1}{24}\left\| {\psi _s
^{\prime \prime \prime }} \right\|_\infty \left| \chi
\right|^3(P\vert \vert Q)}, \,\, \frac{1}{2} \left[
{{f}'(R) - {f}'(r)} \right] V(P\vert \vert Q) \right\}\notag\\
& \leqslant \frac{(R - r)^2}{32}\min \left\{ {\delta _{\Omega _s }
(r,R),\,\,\frac{R - r}{6}\left\| {\psi _s ^{\prime \prime \prime
}} \right\|_\infty },\,\,\frac{R - r}{4} \left[ {{f}'(R) -
{f}'(r)} \right] \right\},\notag
\end{align}

\noindent where
\begin{align}
\label{eq78} E&_{\Omega _s }^\ast (P\vert \vert
Q)\\
& = \begin{cases}
 {\left[ {s(s - 1)} \right]^{ - 1}\sum\limits_{i = 1}^n {(p_i - q_i )\left(
{\frac{p_i + 3q_i }{2(p_i + q_i )}} \right)^s\left( {\frac{p_i + (3 - 2s)q_i
}{p_i + 3q_i }} \right),} } & {s \ne 0,1} \\
 {\sum\limits_{i = 1}^n {(q_i - p_i )\ln \left( {\frac{p_i + 3q_i }{2(p_i +
q_i )}} \right) - \frac{1}{2}\sum\limits_{i = 1}^n {\frac{(p_i - q_i
)^2}{p_i + 3q_i }} ,} } & {s = 0} \\
 {\frac{1}{2}\Delta (P\vert \vert Q) + \frac{1}{2}\sum\limits_{i = 1}^n
{(p_i - q_i )\ln \left( {\frac{p_i + 3q_i }{2(p_i + q_i )}} \right),} } & {s
= 1} \\
\end{cases},\notag
\end{align}

\begin{equation}
\label{eq79} \delta _{^{\Omega _s }} (r,R) = \frac{1}{4}\left[
{\frac{1}{r^3}\left( {\frac{r + 1}{2r}} \right)^{s - 2} -
\frac{1}{R^3}\left( {\frac{R + 1}{2R}} \right)^{s - 2}} \right],
\,\, s \geqslant - 1
\end{equation}

\noindent and
\begin{equation}
\label{eq80} \left\| {\psi _s ^{\prime \prime \prime }}
\right\|_\infty = \frac{(s + 1 + 3r)}{r^2(r + 1)^3}\left( {\frac{r
+ 1}{2r}} \right)^s, \,\, s \geqslant - 1.
\end{equation}
\end{theorem}

\begin{proof} The third order derivative of the
function $\psi _s (x)$ is given by
\begin{equation}
\label{eq81} \psi _s ^{\prime \prime \prime }(x) = - \frac{(s + 1
+ 3x)}{x^2(x + 1)^3}\left( {\frac{x + 1}{2x}} \right)^s, \,\, x
\in (0,\infty )
\end{equation}

This gives
\begin{equation}
\label{eq82} \psi _s ^{\prime \prime \prime }(x) \leqslant 0, \,\,
\forall \mbox{ }s \geqslant - 1.
\end{equation}

From (\ref{eq82}), we can say that the function $\psi ^{\prime
\prime }(x)$ is monotonically decreasing function in $x \in
(0,\infty )$, and hence, for all $x \in [r,R]$, we have
\begin{align}
\label{eq83} \delta _{\Omega _s } (r,R) & = {\psi }''(r) - {\psi
}''(R) \\
& = \frac{1}{4}\left[ {\frac{1}{r^3}\left( {\frac{r + 1}{2r}}
\right)^{s - 2} - \frac{1}{R^3}\left( {\frac{R + 1}{2R}}
\right)^{s - 2}} \right], \,\, s \geqslant - 1.\notag
\end{align}

From (\ref{eq81}), we have
\begin{equation}
\label{eq84} \left| {\psi _s ^{\prime \prime \prime }}
\right|\,^\prime = - \left( {\frac{(s + 1)(s + 2) + 8(s + 1)x +
12x^2}{x^3(x + 1)^4}} \right)\left( {\frac{x + 1}{2x}} \right)^s.
\end{equation}

From (\ref{eq84}), we can say that the function $\left| {\psi _s
^{\prime \prime \prime }} \right|$ is monotonically decreasing
function in $x \in (0,\infty )$ for $s \geqslant - 1$, and hence,
for all $x \in [r,R]$, we have
\begin{equation}
\label{eq85} \left\| {\psi _s ^{\prime \prime \prime }}
\right\|_\infty = \mathop {\sup }\limits_{x \in [r,R]} \left|
{\psi _s ^{\prime \prime \prime }} \right| = \frac{(s + 1 +
3r)}{r^2(r + 1)^3}\left( {\frac{r + 1}{2r}} \right)^s, \,\, s
\geqslant - 1
\end{equation}

Applying Theorem \ref{the33} for the measure (\ref{eq19}) along
with (\ref{eq83}) and (\ref{eq85}) we get the required proof.
\end{proof}

In view of the inequalities (\ref{eq76}) and (\ref{eq77}) we have
some particular cases given in the following corollary.

\begin{corollary} \label{cor44} The following bounds hold:
\begin{align}
\label{eq89} & \left| {\Delta (P\vert \vert Q) - 2\sum\limits_{i =
1}^n {q_i \left( {\frac{p_i - q_i }{p_i + q_i }} \right)^2} }
\right|\\
& \leqslant \min \left\{ {2\left[ {\frac{1}{(r + 1)^3} -
\frac{1}{(R + 1)^3}} \right]\chi ^2(P\vert \vert Q)}
\right.,\notag \\
& \qquad \qquad \frac{4}{(r + 1)^4}\vert \chi \vert ^3(P\vert
\vert Q), \left. {\frac{8(R - r)(R + r + 2)}{(r + 1)^2(R +
1)^2}V(P\vert \vert Q)}
\right\}\notag\\
& \leqslant \min \left\{ {2(R - r)^2\left[ {\frac{1}{(r + 1)^3} -
\frac{1}{(R + 1)^3}} \right]} \right.,\notag\\
& \qquad \qquad \frac{(R - r)^3}{2(r + 1)^4}, \left. {\frac{4(R -
r)^2(R + r + 2)}{(r + 1)^2(R + 1)^2}} \right\}\notag
\end{align}
\begin{align}
\label{eq90} & \left| {F(P\vert \vert Q) - \frac{1}{2}D(Q\vert
\vert P) + \frac{1}{4}\Delta (P\vert \vert Q)} \right|\\
& \leqslant \min \left\{ {\frac{1}{8}\left[ {\frac{1}{r(r + 1)^2}
- \frac{1}{R(R + 1)^2}} \right]\chi ^2(P\vert \vert Q)} \right.,\,
\frac{3r + 1}{12r^2(r + 1)^3}\vert \chi \vert
^3(P\vert \vert Q),\notag\\
& \qquad \qquad \left. {\frac{R - r}{rR}\left[ {L_{ - 1}^{ - 1}
\left( {\frac{r + 1}{r},\frac{R + 1}{R}} \right) - L_{ - 2}^{ - 2}
\left( {\frac{r + 1}{r},\frac{R + 1}{R}} \right)} \right]V(P\vert
\vert Q)} \right\}\notag\\
& \leqslant \min \left\{ {\frac{(R - r)^2}{32}\left[ {\frac{1}{r(r
+ 1)^2} - \frac{1}{R(R + 1)^2}} \right]} \right., \frac{(3r + 1)(R
- r)^3}{96r^2(r + 1)^3},\notag\\
& \qquad \qquad \left. {\frac{(R - r)^2}{2rR}\left[ {L_{ - 1}^{ -
1} \left( {\frac{r + 1}{r},\frac{R + 1}{R}} \right) - L_{ - 2}^{ -
2} \left( {\frac{r + 1}{r},\frac{R + 1}{R}} \right)} \right]}
\right\}\notag
\end{align}

\noindent and
\begin{align}
\label{eq91} & \left| {G(P\vert \vert Q) - \frac{1}{4}\left[ {\chi
^2(Q\vert \vert P) - D(Q\vert \vert P)} \right]} \right|\\
& \leqslant \min \left\{ {\frac{1}{16}\left[ {\frac{1}{r^2(r + 1)}
- \frac{1}{R^2(R + 1)}} \right]\chi ^2(P\vert \vert Q)} \right.,
\frac{3r + 2}{24r^3(r + 1)^2}\vert \chi \vert ^3(P\vert \vert
Q),\notag\\
& \qquad \qquad \left. {\frac{(R - r)}{2rR}\left[ {1 - L_{ - 1}^{
- 1} \left( {\frac{r + 1}{r},\frac{R + 1}{R}} \right)}
\right]V(P\vert \vert Q)} \right\}.\notag \\
& \leqslant \min \left\{ {\frac{(R - r)^2}{64}\left[
{\frac{1}{r^2(r + 1)} - \frac{1}{R^2(R + 1)}} \right]} \right.,
\frac{(3r + 2)(R - r)^2}{192r^3(r + 1)^2},\notag\\
& \qquad \qquad \left. {\frac{(R - r)^2}{4rR}\left[ {1 - L_{ -
1}^{ - 1} \left( {\frac{r + 1}{r},\frac{R + 1}{R}} \right)}
\right]} \right\} .\notag
\end{align}
\end{corollary}

\begin{proof} Letting in (\ref{eq76}), $s = -
1$, $s = 0$ and $s = 1$ we get the inequalities (\ref{eq89}),
(\ref{eq90}) and (\ref{eq91}) respectively.
\end{proof}

In view of the inequalities (\ref{eq77}) we have some particular
cases given in the following corollary.

\begin{corollary} \label{cor45} The following bounds hold:
\begin{align}
\label{eq92} &\left| {\Delta (P\vert \vert Q) - \sum\limits_{i =
1}^n {(p_i + 7q_i )\left( {\frac{p_i - q_i }{p_i + 3q_i }}
\right)^2} } \right|\\
& \leqslant \min \left\{ {\left[ {\frac{1}{(r + 1)^3} -
\frac{1}{(R + 1)^3}} \right]\chi ^2(P\vert \vert Q)}
\right.,\,\frac{1}{(r + 1)^4}\left| \chi \right|^3(P\vert
\vert Q),\notag\\
& \qquad \qquad \left. {\frac{2(R - r)(R + r + 2)}{(r + 1)^2(R +
1)^2}V(P\vert \vert Q)} \right\} .\notag\\
& \leqslant \min \left\{ {\frac{(R - r)^2}{4}\left[ {\frac{1}{(r +
1)^3} - \frac{1}{(R + 1)^3}} \right]} \right., \frac{(R -
r)^3}{8(r + 1)^4}, \notag\\
&\qquad \qquad \left. {\frac{(R - r)^2(R + r + 2)}{(r + 1)^2(R +
1)^2}} \right\} .\notag
\end{align}

\begin{align}
\label{eq93} &\left| {F(P\vert \vert Q) - \sum\limits_{i = 1}^n
{(q_i - p_i )\ln \left( {\frac{p_i - q_i }{p_i + 3q_i }} \right) +
\frac{1}{2}} \sum\limits_{i = 1}^n {\frac{(p_i - q_i )^2}{p_i +
3q_i }} } \right|\\
& \leqslant \min \left\{ {\frac{1}{8}\left[ {\frac{1}{r(r + 1)^2}
- \frac{1}{R(R + 1)^2}} \right]\chi ^2(P\vert \vert Q)} \right.,
\frac{3r + 1}{24r^2(r + 1)^3}\left| \chi \right|^3(P\vert \vert
Q),\notag\\
& \qquad \qquad \left. {\frac{R - r}{2rR}\left[ {L_{ - 1}^{ - 1}
\left( {\frac{r + 1}{r},\frac{R + 1}{R}} \right) - L_{ - 2}^{ - 2}
\left( {\frac{r + 1}{r},\frac{R + 1}{R}} \right)} \right]V(P\vert
\vert Q)}
\right\} .\notag\\
& \leqslant \min \left\{ {\frac{(R - r)^2}{32}\left[ {\frac{1}{r(r
+ 1)^2} - \frac{1}{R(R + 1)^2}} \right]} \right., \frac{(3r + 1)(R
- r)^3}{192r^2(r + 1)^3},\notag \\
& \qquad \qquad \left. {\frac{(R - r)^2}{4rR}\left[ {L_{ - 1}^{ -
1} \left( {\frac{r + 1}{r},\frac{R + 1}{R}} \right) - L_{ - 2}^{ -
2} \left( {\frac{r + 1}{r},\frac{R + 1}{R}} \right)} \right]}
\right\} .\notag
\end{align}

\noindent and
\begin{align}
\label{eq94} & \left| {G(P\vert \vert Q) - \frac{1}{2}\Delta
(P\vert \vert Q) - \frac{1}{2}\sum\limits_{i = 1}^n {(p_i - q_i
)\ln \left( {\frac{p_i - q_i }{p_i + 3q_i }} \right)} } \right|\\
& \leqslant \min \left\{ {\frac{1}{16}\left[ {\frac{1}{r^2(r + 1)}
- \frac{1}{R^2(R + 1)}} \right]\chi ^2(P\vert \vert Q)} \right.,
\frac{3r + 2}{48r^3(r + 1)^2}\left| \chi \right|^3(P\vert \vert
Q),\notag\\
& \qquad \qquad \left. {\frac{(R - r)}{4rR}\left[ {1 - L_{ - 1}^{
- 1} \left( {\frac{r + 1}{r},\frac{R + 1}{R}} \right)}
\right]V(P\vert \vert Q)} \right\}\notag\\
& \leqslant \min \left\{ {\frac{(R - r)^2}{64}\left[
{\frac{1}{r^2(r + 1)} - \frac{1}{R^2(R + 1)}} \right]} \right.,
\frac{(3r + 2)(R - r)^2}{384r^3(r + 1)^2},\notag \\
& \qquad \qquad \left. {\frac{(R - r)^2}{8rR}\left[ {1 - L_{ -
1}^{ - 1} \left( {\frac{r + 1}{r},\frac{R + 1}{R}} \right)}
\right]} \right\}.\notag
\end{align}
\end{corollary}

\begin{proof} Letting in (\ref{eq77}), $s = -
1$, $s = 0$ and $s = 1$ we get the inequalities (\ref{eq92}),
(\ref{eq93}) and (\ref{eq94}) respectively.
\end{proof}

\bigskip
\begin{center}
\textbf{Acknowledgments}
\end{center}

This work has been done during second author's stay with the
\textit{"Mathematics Department, College of Science and
Management, University of Northern British Columbia, Prince George
BC V2N4Z9, Canada"}, for which he is thankful to the above
mentioned university for the support and hospitality.


\begin{thebibliography}{99}
\bibitem{pea} K. PEARSON, On the Criterion that a given system of
eviations from the probable in the case of correlated system of
variables is such that it can be reasonable supposed to have
arisen from random sampling, \textit{Phil. Mag.},
\textbf{50}(1900), 157-172.
\end{thebibliography}
\end{document}